\documentclass{amsart}

\usepackage{amssymb,amsmath}
\usepackage[pdf]{pstricks} 
\usepackage{mathrsfs}
\renewcommand{\mathcal}{\mathscr}
\usepackage{dsfont}
\usepackage{paralist}
\usepackage[all]{xypic}
\xyoption{dvips}
\usepackage{url}
\usepackage[colorlinks,plainpages,backref]{hyperref}

\theoremstyle{definition}
\newtheorem{ntn}{Notation}[section]

\theoremstyle{plain}
\newtheorem{lem}[ntn]{Lemma}
\newtheorem{prp}[ntn]{Proposition}
\newtheorem{thm}[ntn]{Theorem}
\newtheorem{cor}[ntn]{Corollary}

\theoremstyle{remark}

\newtheorem{rmk}[ntn]{Remark}
\newtheorem{exa}[ntn]{Example}

\newcommand{\ld}{\ldots}
\newcommand{\cd}{\cdots}

\newcommand{\pp}{\varphi}
\newcommand{\beq}{\begin{equation}}
\newcommand{\eeq}{\end{equation}}
\newcommand{\ben}{\begin{enumerate}}
\newcommand{\een}{\end{enumerate}}
\newcommand{\y}{\item}

\newcommand{\calC}{\mathcal{C}}
\newcommand{\calF}{\mathcal{F}}
\newcommand{\calK}{\mathcal{K}}
\newcommand{\calO}{\mathcal{O}}
\newcommand{\OO}{\mathcal{O}}
\newcommand{\calR}{\mathcal{R}}
\newcommand{\basering}[1]{\ensuremath{\mathbb{#1}}}
\newcommand{\CC}{\basering{C}}

\newcommand{\eps}{\varepsilon}
\newcommand{\fg}{\mathfrak{g}}
\newcommand{\ideal}[1]{{\left\langle #1\right\rangle}}
\newcommand{\into}{\hookrightarrow}
\newcommand{\onto}{\twoheadrightarrow}
\newcommand{\mm}{\mathfrak{m}}
\newcommand{\mmm}{m^\mu_\mu}
\newcommand{\p}{\partial}
\newcommand{\rel}{\mathrm{rel}}
\newcommand{\xymat}{\SelectTips{cm}{}\xymatrix}
\newcommand{\rdots}{\mathinner{%
  \mkern1mu\raise1pt\hbox{.}%
  \mkern2mu\raise4pt\hbox{.}%
  \mkern2mu\raise7pt\vbox{\kern7pt\hbox{.}}\mkern1mu}} 

\DeclareMathOperator{\Ann}{Ann}
\DeclareMathOperator{\depth}{depth}
\DeclareMathOperator{\diag}{diag}
\DeclareMathOperator{\Der}{Der}
\DeclareMathOperator{\Hom}{Hom}
\DeclareMathOperator{\GL}{GL}
\DeclareMathOperator{\Sing}{Sing}

\numberwithin{equation}{section}

\begin{document}

\title{Adjoint divisors and free divisors}

\author{David Mond}
\address{
D. Mond\\
Mathematics Institute\\
University of Warwick\\
Coventry CV47AL\\
England
}
\email{D.M.Q.Mond@warwick.ac.uk}

\author{Mathias Schulze}
\address{
M. Schulze\\
Department of Mathematics\\
University of Kaiserslautern\\
67663 Kaiserslautern\\
Germany}
\email{mschulze@mathematik.uni-kl.de}
\thanks{The research leading to these results has received funding from the People Programme (Marie Curie Actions) of the European Union's Seventh Framework Programme (FP7/2007-2013) under REA grant agreement n\textsuperscript{o} PCIG12-GA-2012-334355.}

\date{\today}


\begin{abstract}
We describe two situations where adding the adjoint divisor to a divisor $D$ with smooth normalization yields a free divisor. 
Both also involve stability or versality.
In the first, $D$ is the image of a corank $1$ stable map-germ $(\CC^n,0)\to(\CC^{n+1},0)$, and is not free.
In the second, $D$ is the discriminant of a versal deformation of a weighted homogeneous function with isolated critical point (subject to certain numerical conditions on the weights). 
Here $D$ itself is already free. 

We also prove an elementary result, inspired by these first two, from which we obtain a plethora of new examples of free divisors.
The presented results seem to scratch the surface of a more general phenomenon that is still to be revealed. 
\end{abstract}

\subjclass{32S25, 17B66, 32Q26, 32S30, 14M17}
\keywords{free divisor, adjoint divisor, stable map, versal deformation, isolated singularity, prehomogeneous vector space}

\maketitle
\tableofcontents

\section{Introduction}

Let $M$ be an $n$-dimensional complex analytic manifold and $D$ a hypersurface in $M$.
The $\OO_M$-module $\Der(-\log D)$ of \emph{logarithmic vector fields} along $D$ consists of all vector fields on $M$ tangent to $D$ at all smooth points of $D$.
If this module is locally free, $D$ is called a \emph{free divisor}.
This terminology was introduced by Kyoji Saito in \cite{Sai80b}. 
As freeness is evidently a local condition, so we may pass to germs of analytic spaces $D\subset X:=(\CC^n,0)$, pick coordinates $x_1,\dots,x_n$ on $X$ and a defining equation $h\in\OO_X=\CC\{x_1,\dots,x_n\}$ for $D$.

The module $\Der(-\log D)$ is an infinite-dimensional Lie sub-algebra of $\Der_X$, to be more precise, a Lie algebroid.
Thus free divisors bring together commutative algebra, Lie theory and the theory of ${\mathscr D}$-modules, see \cite{CN05}. 
Freeness has been used by Jim Damon and the first author to give an algebraic method for computing the vanishing homology of sections of discriminants and other free divisors, see \cite{DM91} and \cite{Dam96}. 
More recently the idea of adding a divisor to another in order to make the union free has been used by Damon and Brian Pike as a means of extending this technique to deal with sections of non-free divisors, see \cite{DP11} and \cite{DP12}. 

Saito formulated the following elementary freeness test, called \emph{Saito's criterion} (see \cite[Thm.~1.8.(ii)]{Sai80b}):
If the determinant of the so-called \emph{Saito matrix} $(\delta_i(x_j))$ generates the defining ideal $\ideal{h}$ for some $\delta_1,\ldots,\delta_n\in\Der(-\log D)$, then $D$ is free and $\delta_1,\ldots,\delta_n$ is a basis of $\Der(-\log D)$.
While any smooth hypersurface is free, singular free divisors are in fact highly singular:
Let $\Sing D$ be the singular locus of $D$ with structure defined by the Jacobian ideal of $D$.
By the theorem of Aleksandrov--Terao (see \cite[\S1 Thm.]{Ale88} or \cite[Prop.~2.4]{Ter80a}), freeness of $D$ is equivalent to $\Sing D$ being a Cohen--Macaulay space of (pure) codimension $1$ in $D$.

The simplest example of a free divisor, whose importance in algebraic geometry is well-known, is the normal crossing divisor $D$ defined by $h:=x_1\cdots x_n$; here, due to Saito's criterion, $\Der(-\log D)$ is freely generated by the vector fields $x_1\frac\p{\p x_1},\ldots, x_n\frac\p{\p x_n}$.
In general free divisors are rather uncommmon: given $n$ vector fields $\delta_1,\dots,\delta_n\in\Der_X$, let $h$ be the determinant of their Saito matrix, and suppose that $h$ is reduced. 
Then $h$ defines a free divisor if and only if the $\OO_{\CC^n}$-submodule of $\Der_M$ generated by the $\delta_j$ is a Lie algebra, see \cite[Lem.~1.9]{Sai80b}. 
Thus to generate examples, special techniques are called for.
Non-trivial examples of free divisors first appeared as discriminants and bifurcation sets in the base of versal deformations of isolated hypersurface singularities, see \cite{Sai80a}, \cite{Ter83}, \cite{Loo84}, \cite{Bru85}, \cite{vSt95}, \cite{Dam98}, \cite{BEGvB09}. 
Here freeness follows essentially from the fact that $\Der(-\log D)$ is the kernel of the \emph{Kodaira--Spencer map} from the module of vector fields on the base to the relative $T^1$ of the deformation.
\smallskip

In this paper, we construct new examples of free divisors by a quite different procedure. Recall that 
we denote by $X$ the germ $(\CC^n,0)$. 
Let now $f\colon X\to(\CC^{n+1},0)=:T$ be a finite and generically $1$-to-$1$ holomorphic map germ. 
In particular, $X$ is a normalization of the reduced image $D$ of $f$.
Denote by $\calF_i$ the $i$th Fitting ideal of $\OO_X$ considered as $\OO_T$-module.
Mond and Pellikaan~\cite[Props.~3.1, 3.4, 3.5]{MP89} showed that $D$ is defined by $\calF_0$, $\calF_1$ is perfect ideal of height $2$ restricting to the \emph{conductor ideal} $\calC_D:=\Ann_{\OO_D}(\OO_X/\OO_D)$ which in turn is a principal ideal of $\OO_X$.
In particular, the reduced singular locus $\Sigma$ of $D$ is the closure of the set of double points of $f$ and $\calC_D$ is radical provided $D$ is normal crossing in codimension $1$.
We call any member of $\calF_1$ whose pull-back under $f$ generates $\calC_D$ an \emph{adjoint equation}, and its zero locus $A$ an \emph{adjoint divisor}. 
Set-theoretically, this implies that $\Sigma=A\cap D$.
Ragni Piene~\cite[\S3]{Pie79} showed that an adjoint equation is given by the quotient
\[
\frac{\p h/\p x_j}{\p(f_1,\dots,\widehat {f_j},\ld, f_{n+1})/\p(x_1,\ld,x_n)}.
\]
For example, if $f(x_1, x_2)=(x_1, x_2^2, x_1x_2)$ is the parameterisation of the Whitney umbrella, with image  $D=\{t_3^2-t_1^2t_2=0\}$, then this recipe gives $t_1$ as adjoint equation. 
Figure~\ref{22} below shows $D+A$ in this example.

We show
\begin{thm}\label{59}
Let $D$ be the image of a stable map germ $(\CC^n,0)\to(\CC^{n+1},0)$ of corank $1$, and let $A$ be an adjoint divisor for $D$. 
Then $D+A$ is a free divisor.
\end{thm}

However, we show by an example that $D+A$ is not free when $D$ is the image of a stable germs of corank $\geq 2$ (see Example~\ref{17}.\eqref{17c}).
We recall the standard normal forms for stable map germs of corank $1$ in \S\ref{images}.
\smallskip

In \S\ref{discr}, we prove the following analogous result for the discriminants of certain weighted homogeneous isolated function singularities. Recall that the discriminant has smooth normalization, so that the preceding definitions apply. 

\begin{thm}\label{60} 
Let $f\colon(\CC^n,0)\to(\CC,0)$ be a weighted homogeneous polynomial of degree $d$ with isolated critical point and Milnor number $\mu$.
Let $d_1\geq d_2\geq\cdots\geq d_\mu$ denote the degrees of the members of a weighted homogeneous $\CC$-basis of the Jacobian algebra $\OO_{\CC^n,0}/\ideal{\frac{\p f}{\p x_1},\dots,\frac{\p f}{\p x_n}}$.
Assume that $d-d_1+2d_i\neq 0\neq d-d_i$ for $i=1,\ldots,\mu$.

Let $D$ be the discriminant in the base space of an $\calR_e$-versal deformation of $f$ and let $A$ be an adjoint divisor for $D$. 
Then $D+A$ is a free divisor. 
\end{thm}

Theorem~\ref{60} evidently applies to the simple singularities, since for these $d_1<d$. 
It also applies in many other cases.
For example, it is easily checked that the hypotheses on the weights hold for plane curve singularities of the form $x^p+y^q$ with $p$ and $q$ coprime. 
We do not know whether the conclusion continues to holds without the hypotheses on the weights. 

The adjoint $A$ of the discriminant is closely related to the bifurcation set; this is discussed at the end of Section \ref{discr}.

We remark that an adjoint can be defined verbatim in case $X$ is merely Gorenstein rather than smooth.
In \cite[\S6]{GMS12} the techniques used here are applied to construct new free divisors from certain Gorenstein varieties lying canonically over the discriminants of Coxeter groups. 

In Theorem~\ref{60}, $D$ is already a free divisor as remarked above. 
In contrast, in Theorem \ref{59}, $D$ itself is not free: the argument with the Kodaira Spencer map referred to above shows that $\Der(-\log D)$ has depth $n$ rather than $n+1$. 
So by adding $A$ we are making a non-free divisor free. 
Something similar was already done by Jim Damon, for the same divisor $D$, in \cite[Ex.~8.4]{Dam98}; Damon showed that after the addition of a certain divisor $E$ (not an adjoint divisor for $D$) with two irreducible components, $D+E$ is the discriminant of a $\calK_V$-versal deformation of a non-linear section of another free divisor $V$. 
Freeness of $D+E$ followed by his general theorem on $\calK_V$-versal discriminants. 
It seems that our divisor $D+A$ does not arise as a discriminant using Damon's procedure.

\smallskip
 
A crucial step in the proofs of both of Theorems~\ref{59} and \ref{60} is the following fact (see Propositions~\ref{104} and \ref{13}). 

\begin{prp}\label{cycl}
In the situations of Theorems~\ref{59} and \ref{60}, $\calF_1$ is cyclic as module over the Lie algebra $\Der(-\log D)$, and generated by an adjoint equation $a\in\OO_{\CC^{n+1},0}$ for $D$.
\end{prp}

Indeed, Theorem~\ref{60} follows almost trivially from this (see Proposition~\ref{106} below). 
We cannot see how to deduce Theorem~\ref{59} in an equally transparent way.
Unfortunately our proof of Proposition~\ref{cycl} is combinatorial and not very revealing, something we hope to remedy in future work.

From Proposition~\ref{cycl} it follows that the adjoint is unique up to isomorphism preserving $D$ (see Corollaries~\ref{111} and \ref{115}).

In the process of proving Theorem~\ref{60}, we note an easy argument which shows that the preimage of the adjoint divisor in the normalization of the discriminant is itself a free divisor. 
This is our Theorem~\ref{12}. 

\smallskip

Motivated by Theorems~\ref{59} and \ref{60}, we describe in \S\ref{lfd} a general procedure which constructs, from a triple consisting of a free divisor $D$ with $k$ irreducible components and a free divisor in $(\CC^k,0)$ containing the coordinate hyperplanes, a new free divisor containing $D$.
By this means we are able to construct a surprisingly large number of new examples of free divisors. 

\subsection*{Notation}

We shall denote by $\calF_\ell^R(M)$ the $\ell$th \emph{Fitting ideal} of the $R$-module $M$.
For any presentation
\[
\xymat{
R^m\ar[r]^-A & R^k\ar[r] & M \ar[r] & 0
}
\]
it is generated by all $(k-\ell)$-minors of $A$ and defines the locus where $M$ requires more than $\ell$ generators.

For any analytic space germ $X$, we denote by $\Theta_X:=\Hom_{\calO_X}(\Omega_X^1,\calO_X)$ the $\calO_X$-module of vector fields on $X$.
The $\calO_X$-module of vector fields along an analytic map germ $f\colon X\to Y$ is defined by $\Theta(f):=f^*\Theta_Y$.
For $X$ and $Y$ smooth, we shall use the standard operators $tf:\Theta_X\to \Theta(f)$ and $\omega f:\Theta_Y\to\Theta(f)$ of singularity theory defined by $tf(\xi):=Tf\circ\xi$ and $\omega f(\eta)=\eta\circ f$.
For background in singularity theory we suggest the survey paper \cite{Wal81} of C.T.C.~Wall.

Throughout the paper, all the hypersurfaces we consider will be assumed reduced, without further mention.

\subsection*{Acknowledgments} 

We are grateful to Eleonore Faber for pointing out a gap in an earlier version of Theorem~\ref{19}, and to the referee for a number of valuable suggestions.

\section{Images of stable maps}\label{images}

Let $f\colon X:=(\CC^n,0)\to(\CC^{n+1},0)=:T$ be a finite and generically one-to-one map-germ with image $D$. 
Note that $X=\bar D$ is a normalization of $D$.
By \cite[Prop.~2.5]{MP89}, the $\calO_T$-module $\calO_X$ has a free resolution of the form 
\beq\label{23}
\xymat{
0\ar[r] & \calO_T^k\ar[r]^-\lambda & \calO_T^k\ar[r]^-\alpha 
& \calO_X\ar[r] & 0,
}
\eeq
in which the matrix $\lambda$ can be chosen symmetric (we shall recall the proof below). 
For $1\le i,j\le k$, we denote by $m^i_j$ the minor obtained from $\lambda$ by deleting the $i$th row and the $j$th column.
The map $\alpha$ sends the $i$th basis vector $e_i$ to $g_i\in\calO_X$, where $g=g_1,\dots,g_k$ generates $\calO_X$ over $\calO_T$. 
It will be convenient to assume, after reordering the $g_i$, that $g_k=1$. 
This leads to a free presentation
\[
\xymat{
\calO_T^k\ar[r]^-{\lambda^{\hat k}} & \calO_T^{k-1}\ar[r]^-\alpha 
& \OO_X/\OO_D\ar[r] & 0,
}
\]
where $\lambda^{\hat k}$ is obtained from $\lambda$ by deleting the $k$th row (corresponding to the generator $g_k=1$ of $\OO_X$).
By a theorem of Buchsbaum and Eisenbud~\cite{BE77}, this shows that 
\begin{equation}\label{126}
\calF_1':=\calF_0^{\OO_T}(\OO_X/\OO_D)=\Ann_{\OO_T}(\OO_X/\OO_D)=\calC_D\OO_T.
\end{equation}
As $\calF_\ell^{\OO_T}(\OO_X)$ defines the locus where $\OO_X$ requires more than $\ell$ $\OO_T$-generators, $\det(\lambda)$ is an equation for $D$.
By the hypothesis that $f$ is generically one-to-one, it is a reduced equation for $D$ (see \cite[Prop.~3.1]{MP89}).

By Cramer's rule one finds that in $\calO_X$, $g_im^j_s=\pm g_jm^i_s$ for $1\leq i,j,s\leq k$ (see~\cite[Lem.~3.3]{MP89}); invoking the symmetry of $\lambda$, this gives
\begin{equation}\label{124}
g_im^k_j=\pm g_jm^k_i,\quad 1\leq i,j\leq k.
\end{equation} 
Combining this with the structure equations $g_ig_j=\sum_{\ell=1}^k\alpha^\ell_{i,j}g_\ell$ for $\OO_X$ as $\OO_T$-algebra, one shows that all of the $m^i_j$ lie in $\Ann_{\OO_T}(\OO_X/\OO_D)$; then from \eqref{126}, one deduces that (see \cite[Thm. 3.4]{MP89})

\begin{lem}\label{109}
$\calF_1=\calF_1'=\ideal{m^\mu_j\mid j=1,\dots,\mu}$ is a determinantal ideal.
\end{lem}

Since $\calF_1$ therefore gives $\Sing D$ a Cohen-Macaulay structure, this structure is reduced if generically reduced. 
If $f$ is stable, then at most points of $\Sing D$, $D$ consists of two smooth irreducible components meeting transversely, from which generic reducedness follows. 

As $g_k=1$, it follows from \eqref{124} that $m^k_i=\pm g_im^k_k$ and hence

\begin{lem}\label{110} 
$\calF_1\OO_X=\ideal{m^k_k}_{\calO_X}$ is a principal ideal.
\end{lem}

From this, we deduce

\begin{lem}\label{107}
Any adjoint divisor $A$ for $D$ is of the form 
\[
A=V\bigl(m^k_k+\sum_{j=1}^{k-1}c_jm^k_j\bigr),\quad c_j\in\calO_T.
\]
\end{lem}

Mather~\cite{Mat69b} showed how to construct normal forms for stable map germs: one begins with a germ $f\colon(\CC^k,0)\to(\CC^\ell,0)$ whose components lie in $\mm_{\CC^k,0}^2$, and unfolds it by adding to $f$ terms of the form $u_ig_i$, where the $u_i$ are unfolding parameters and the $g_i$ form a basis for $\mm_{\CC^k,0}\Theta(f)/(tf(\Theta_{\CC^k,0})+\mm_{\CC^\ell,0}\Theta(f))$. 
Applying this construction to $f(x)=(x^{k},0)$, one obtains the stable corank-$1$ map germ $f_k:(\CC^{2k-2},0)\to (\CC^{2k-1},0)$ given by
\beq\label{nf}
f_k(u,v,x):=(u,v,x^{k}+u_1x^{k-2}+\cdots+u_{k-2}x, v_1x^{k-1}+\cdots+v_{k-1}x)=(u,v,w),
\eeq
where we abbreviate $u:=u_1,\ldots,u_{k-2}$, $v:=v_1,\ldots,v_{k-1}$, $w:=w_1, w_2$.

\begin{exa}\label{17}\
\begin{asparaenum}

\item\label{17a} When $k=2$ in \eqref{nf}, $f_2(v,x)=(v,x^2,vx)$ parameterizes the Whitney umbrella in $(\CC^3,0)$, whose equation is $w_2^2-v^2w_1=0$.
With respect to the  basis $g=x,1$ of $\OO_S$ over $\OO_T$, one has the symmetric presentation matrix
\[
\lambda=
\begin{pmatrix}
v&-w_2\\
-w_2&vw_1
\end{pmatrix},
\]
and this gives $a=v$ as equation of an adjoint $A$ (see Figure~\ref{22}). 

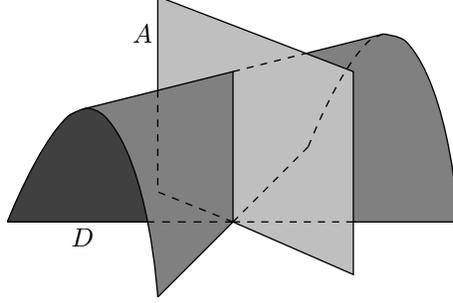
\begin{figure}[ht]
\caption{Whitney umbrella with adjoint}\label{22}
\begin{center}
\begin{pspicture}(6,4)
\rput(3,2){
\psset{dash=1mm 1mm}
\pscustom[fillcolor=gray,fillstyle=solid,linestyle=none]{
\pscurve(2,1.5)(2.3,1.4)(3,-1)
\psline(3,-1)(0,-1)(0,1)(2,1.5)}
\pscustom[fillcolor=lightgray,fillstyle=solid,linestyle=none]{
\psline(-1,2)(1.6,1)(1.6,-1.7)(0,-1)(-1,-0.6)
}
\psset{linewidth=0pt}
\pscustom[fillcolor=darkgray,fillstyle=solid,linestyle=none]{
\pscurve(-2,0.5)(-2.2,0.4)(-3,-1)
\psline(-3,-1)(0,-1)(0,1)}
\psset{linewidth=0.5pt}
\pscustom[fillcolor=gray,fillstyle=solid]{
\psline(-1,-2)(0,-1)(0,1)(-2,0.5)
\pscurve(-2,0.5)(-1.7,0.4)(-1,-2)}
\pscurve(-2,0.5)(-2.2,0.4)(-3,-1)
\psline(-3,-1)(-1.15,-1)
\psset{linestyle=dashed}
\psline(-1.15,-1)(0,-1)
\psset{linestyle=dashed}
\psline(0,-1)(1.6,-1)
\psset{linestyle=solid}
\psline(1.6,-1)(3,-1)
\pscurve(2,1.5)(2.3,1.4)(3,-1)
\psset{linestyle=dashed}
\pscurve(2,1.5)(1.8,1.4)(1,0)
\psline(1,0)(0,-1)
\psset{linestyle=solid}
\psline(2,1.5)(1,1.25)
\psset{linestyle=dashed}
\psline(1,1.25)(0,1)
\psset{linestyle=solid}
\psline(-1,0.75)(-1,2)(1.6,1)(1.6,-1.7)(0,-1)
\psset{linestyle=dashed}
\psline(0,-1)(-1,-0.6)(-1,0.75)
\rput(-2,-1.2){$D$}
\rput(-1.2,1.5){$A$}
}
\end{pspicture}
\end{center}
\end{figure}

One calculates that $\Der(-\log D)$ and $\Der(-\log (D+A))$ are generated, respectively, by the vector fields whose coefficients are displayed as the columns of the matrices 
\[
\begin{pmatrix}
v&v&0&w_2\\
2w_1&0&2w_2&0\\
0&w_2&v^2&vw_1
\end{pmatrix}
\quad\text{and}\quad
\begin{pmatrix}
v&v&0\\
2w_1&0&2w_2\\
2w_2&w_2&w^2
\end{pmatrix}.
\]
Note that here the basis of the free module $\Der(-\log (D+A))$ is included in a (minimal) list of generators of $\Der(-\log D)$. 
As we will see, this is always the case for the germs $f_k$ described above.

\item\label{17b} Let $D_0$ be the image of the stable corank-$2$ map germ 
\[
(u_1,u_2,u_3,u_4,x,y)\mapsto (u_1,u_2,u_3,u_4,x^2+u_1y,xy+u_2x+u_3y,y^2+u_4x).
\]
One calculates that $D_0+A$ is \emph{not} free.

\item\label{17c} Every stable map germ of corank $k>2$ is adjacent to the germ considered in the preceding example. 
That is, there are points on the image $D$ where $D$ is isomorphic to the product of the divisor $D_0$ we have just considered in \eqref{17b} and a smooth factor. 
It follows that that $D+A$ also is not free.

\item\label{17e} If $D$ is the image of an unstable corank $1$ germ then in general $D+A$ is not free.

\item\label{17d}
A multi-germ of immersions $(\CC^n,\{p_1,\dots,p_k\})\to(\CC^{n+1},0)$ is stable if and only if it is a normal crossing. 
There are strata of such normal crossing points, of different multiplicities, on the image $D$ of a stable map germ such as \eqref{nf}. 
It is easy to show that adding an adjoint divisor to $D$ gives a free divisor at such normal crossing points. 
The normal crossing divisor $D=V(y_1\cdots y_k)\subset(\CC^k,0)$ is normalized by separating its irreducible components.
The equation $a:=\sum_{j=1}^ky_1\cdots\hat{y_j}\cdots y_k$ restricts to $y_1\cdots\hat y_j\cdots y_k$ on the component $V(y_j)$, and so generates the conductor there, and thus $A:=V(a)$ is an adjoint for $D$.
The Euler vector field $\chi=\sum_{j=1}^ky_j\p_{y_j}$, together with the $n-1$ vector fields $\delta_j:=y_j^2\p_{y_j}-y_{j+1}^2\p_{y_{j+1}}$, $j=1,\ldots, k-1$, all lie in $\Der(-\log (D+A))$. 
An application of Saito's criterion yields freeness of $D+A$ and shows that $\chi, \delta_1,\ldots,\delta_{k-1}$ form a basis for $\Der(-\log (D+A))$. 

\end{asparaenum}
\end{exa}

Our proof of Theorem~\ref{59} is based on Saito's criterion.
By Mather's construction, we are concerned with the map $f:=f_k$ of \eqref{nf} where now $n=2k-2$.
Using an explicit list of generators of $\Der(-\log D)$ constructed by Houston and Littlestone in \cite{HL09}, and testing them on the equation $m^k_k$ of $A$, we find a collection of vector fields $\xi_1,\dots,\xi_{2k-1}$ in $\Der(-\log D)$ which are in $\Der(-\log A)$ ``to first order'', in the sense that for $j=1,\dots,2k-1$, we have $\xi_j\cdot m^k_k\in\ideal{m^k_k}+\mm_T\calF_1$. 
Note that $\calF_1$ is intrinsic to $D$, and therefore invariant under any infinitesimal automorphism of $D$, so that necessarily $\xi_j\cdot m^k_k\in\calF_1$. 
In the process of testing, we show that the map $\Der(-\log D)\to \calF_1$ sending $\xi$ to $\xi\cdot m^k_k$ is surjective. 
Using this, we can then adjust the $\xi_j$, without altering their linear part, so that now $\xi_j\cdot m^k_k\in\ideal{m^k_k}$ for $j=1,\dots, 2k-1$. 
As a consequence, the determinant of their Saito matrix must be divisible by the equation of $D+A$.
This determinant contains a distinguished monomial also present in the equation of $D+A$, so the quotient of the determinant by the equation of $D+A$ is a unit, the determinant is a reduced equation for $D+A$, and $D+A$ is a free divisor, by Saito's criterion.

To begin this process, we need more detailed information about the matrix $\lambda$ of \eqref{23}. 
We use a trick from \cite[\textsection 2]{MP89}:
embed $X$ as $X\times\{0\}$ into $X\times(\CC,0):=S$, and let the additional variable in $S$ be denoted by
$t$.  
Extend $f\colon X\to T$ to a map $F\colon S\to T$ by adding $t$ to the last component of $f$. 
Applying this procedure to the map $f=f_k$ of \eqref{nf}, gives
\begin{gather}\label{nF}
F(u,v,x,t):=(u,v,x^{k}+u_1x^{k-2}+\cdots+u_{k-2}x,v_1x^{k-1}+\cdots+v_{k-1}x+t)=(u,v,w).
\end{gather}
It is clear that $\calO_S/\mm_T\calO_S$ is generated over $\CC$ by the classes of $1,x,\dots,x^{k-1}$, from which it follows by \cite[Cor.~2, p.~137]{Gun74} that $\calO_S$ is a free $\calO_T$-module on the basis 
\beq\label{105}
g_i:=x^{k-i},\quad i=1,\dots,k.
\eeq
Consider the diagram 
\beq\label{tricky}
\xymatrix{0\ar[r]&\OO_S\ar[r]^t&\OO_S\ar[r]&\OO_X\ar[r]&0\\
0\ar[r]&\OO_T^k\ar[r]_{[t]^g_g}\ar[u]_{\pp_g}&\OO_T^k\ar[r]\ar[u]_{\pp_g}&\OO_X\ar@{=}[u]\ar[r]&0}
\eeq
in which $\pp_g$ is the $\OO_T$-isomorphism sending $(c_1,\ld,c_k)\in\OO_T^k$ to $\sum_{j=1}^kc_jg_j\in\OO_S$, and 
where now $[t]^g_g$ denotes the matrix of multiplication by $t$ with respect to the basis $g$ of $\OO_S$ as $\OO_T$-module. The lower row is thus a presentation of $\OO_X$ as $\OO_T$-module. 
This can be improved by a change of basis on the source of $[t]^g_g$, as follows.

Since $\calO_S$ is Gorenstein, $\Hom_{\calO_T}(\calO_S,\calO_T)\cong\calO_S$ as $\calO_S$-module.
Let $\Phi$ be an $\calO_S$-generator of $\Hom_{\calO_T}(\calO_S,\calO_T)$. 
It induces a symmetric perfect pairing 
\[
\ideal{\cdot,\cdot}:\calO_S\times\calO_S\to\calO_T,\quad\ideal{a,b}=\Phi(ab)
\]
with respect to which multiplication by $t$ is self-adjoint.
We refer to this pairing as the Gorenstein pairing.
Now choose a basis $\check g=\check g_1,\dots,\check g_k$ for $\calO_S$ over $\calO_T$ dual to $g$ with respect to $\ideal{\cdot,\cdot}$; that is, such that $\ideal{g_i,\check g_j}=\delta_{i,j}$. 
Then the $(i,j)$th entry of $[t]^{\check g}_g$ equals $\ideal{t\check g_j,\check g_i}$, and so redefining
\beq\label{106}
\lambda=(\lambda^i_j)=(\lambda_1,\dots,\lambda_k):=[t]^{\check g}_g
\eeq
yields a symmetric presentation matrix.

\begin{lem}\label{28}
With an appropriate choice of generator $\Phi$ of $\Hom_{\calO_T}(\calO_S,\calO_T)$, we have
\beq\label{41}
\lambda\equiv
\begin{pmatrix}
-v_1&-v_2&-v_3&\cdots&-v_{k-1}&w_2\\
-v_2&-v_3&&\rdots&w_2&0\\
-v_3&&&\rdots&\rdots&\vdots\\
\vdots&\rdots&\rdots&&&\vdots\\
-v_{k-1}&w_2&\rdots&&&\vdots\\
w_2&0&\cdots&&\cdots&0
\end{pmatrix}\mod\ideal{u,w_1}.
\eeq
\end{lem}

\begin{proof} 
Let $b\in\OO_S$ map to the socle of the $0$-dimensional Gorenstein ring $\calO_S/\mm_T\calO_S$.
This means that
\[
\mm_Sb\subset\mm_T\calO_S.
\]
Thus, for any $a\in\mm_S$ and $\Psi\in\Hom_{\calO_T}(\calO_S,\calO_T)$, we have 
\[
(a\Psi)(b)=\Psi(ab)=\Psi(\sum_ic_ib_i)=\sum_ic_i\Psi(b_i)\in\mm_T
\]
where $c_i\in\mm_T$ and $b_i\in\OO_S$.
In other words, any $\Phi\in\mm_S\Hom_{\calO_T}(\calO_S,\calO_T)$ maps $b$ to $\mm_T$.
Conversely, for any $\Phi\in\Hom_{\calO_T}(\calO_S,\calO_T)$ with $\Phi(b)\in\OO_T^*$, its class $\bar\Phi$ in $\Hom_{\calO_T}(\calO_S,\calO_T)/\mm_S\Hom_{\calO_T}(\calO_S,\calO_T)$ is non-zero.
As $\Hom_{\calO_T}(\calO_S,\calO_T)\cong\OO_S$, this latter space is $1$-dimensional over $\CC$ which shows that $\bar\Phi$ is a generator.
By Nakayama's lemma, $\Phi$ is then an $\OO_S$-basis of $\Hom_{\calO_T}(\calO_S,\calO_T)$.
In particular, we may take as $\Phi(h)$, for $h\in\calO_S$, the coefficient of $b=x^{k-1}$ in the representation of $h$ in the basis \eqref{105}; then $\Phi(b)=1$.

In the following, we implicitly compute modulo $\ideal{u,w_1}$.
Using the relation
\beq\label{45}
w_1=x^k+\sum_{i=1}^{k-2}u_ix^{k-1-i}
\eeq
from \eqref{nF} we compute
\[
\check g_1=1,\ \check g_j=\left(w_1-\sum_{i=j-1}^{k-2}u_ix^{k-1-i}\right)\Big/x^{k+1-j}
=x^{j-1}+\sum_{i=1}^{j-2}u_ix^{j-i-2},\ j=2,\dots,k.
\]
Note that 
\beq\label{30}
\check g_2=x\check g_1,\quad\check g_j=x\check g_{j-1}+u_{j-2}\check g_1,\ j=3,\dots,k.
\eeq
Now let us calculate the columns $\lambda_1,\dots,\lambda_k$ of the matrix \eqref{106}.
Using $\check g_1=1$ and the relation $t+\sum_{i=1}^{k-1}v_ig_i=w_2=w_2g_k$ from \eqref{nF}, we first compute
\[
\lambda_1= [t\check g_1]_g=(\ideal{t,\check g_j})_j=\left(\ideal{w_2g_k-\sum_{i=1}^{k-1}v_ig_i,\ \check g_j}\right)_j=(-v_1,\dots,-v_{k-1}, w_2)^t.
\]
By \eqref{30}, each of the remaining columns $\lambda_j$ is obtained by multiplying by $x$ the vector represented by its predecessor, and, for $j\ge3$, adding $u_{j-2}\lambda_1$. 
Thus, 
\beq\label{44}
\lambda_2=[x]^g_g\lambda_1,\quad\lambda_j=[x]^g_g\lambda_{j-1}+u_{j-2}\lambda_1,\ j=3,\dots,k.
\eeq 
Using \eqref{45} again, observe that
\beq\label{47}
[x]^g_g=
\begin{pmatrix}
w_1&1&0&\cdots&0\\
0&0&1&\ddots&\vdots\\
-u_1&\vdots&\ddots&\ddots&0\\
\vdots&\vdots&&\ddots&1\\
-u_{k-2}&0&\cdots&\cdots&0
\end{pmatrix}
\equiv
\begin{pmatrix}
0&1&0&\cdots&0\\
0&0&1&\ddots&\vdots\\
\vdots&\vdots&\ddots&\ddots&0\\
\vdots&\vdots&&\ddots&1\\
0&0&\cdots&\cdots&0
\end{pmatrix}
\mod\ideal{u,w_1}.
\eeq
The result follows.
\end{proof}

\begin{cor}\label{113}
The reduced equation $h$ of the image $D$ of the map $f_k$ of \eqref{nf} contains the monomial $w_2^k$ with coefficient $\pm 1$. 
The minor $m^k_k$ contains the monomial $w_2^{k-2}v_1$ with coefficient $\pm 1$.
\end{cor}

\begin{proof}
The determinant of the matrix $\lambda$ of \eqref{41} is a reduced equation for the image of $f$ (see \cite[Prop.~3.1]{MP89}).
Both statements then follow from Lemma~\ref{28}.
\end{proof}

\begin{exa}
For the stable map-germs
\[
f_3(u_1,v_1,v_2,x)= (u_1,v_1,v_2,x^3+u_1x,v_1x^2+v_2x)
\]
and
\[
f_4(u_1,u_2,v_1,v_2,v_3,x)=(u_1,u_2,v_1,v_2, v_3,x^4+u_1x^2+u_2x, v_1x^3+v_2x^2+v_3x)
\]
of \eqref{nf}, the matrix $\lambda$ is equal to
\[
\begin{pmatrix}
-v_1&-v_2&w_2\\
-v_2&w_2+u_1v_1&-v_1w_1\\
w_2&-v_1w_1&v_2w_1-u_1w_2
\end{pmatrix}
\]
and to
\[
\begin{pmatrix}
-v_1&-v_2&-v_3&w_2\\
-v_2&u_1v_1-v_3&w_2+u_2v_1&-v_1w_1\\
-v_3&w_2+u_2v_1&u_2v_2-u_1v_3-v_1w_2&-v_2w_1+u_1w_2\\
w_2&-v_1w_2&-v_2w_1+u_1w_2&-v_3w_1+u_2w_2
\end{pmatrix}
\]
respectively.
\end{exa}

\begin{proof}[Proof of Theorem~\ref{59}]
In \cite[Thms.~3.1-3.3]{HL09}, Houston and Littlestone give an explicit list of generators for $\Der(-\log D)$. 
Their proof that these generators lie in $\Der(-\log D)$ simply exhibits, for each member $\xi$ of the list, a lift $\eta\in\Theta_X$ in the sense that $tf(\eta)=\omega f(\xi)$.
The Houston-Littlestone list consists of the Euler field $\xi_e$ and three families $\xi^i_j$, $1\leq i\leq 3$, $1\leq j\leq k-1$.
Denote by $\bar\xi^i_j$ the linear part of $\xi_j^i$. 
After dividing by $1$, $k$, $k$, $k$ and $k^2$ respectively, these linear parts are
\begin{gather*}
\bar\xi_e=
\sum_{i=1}^{k-2}(i+1)u_i\p_{u_i}+\sum_{i=1}^{k-1}iv_i\p_{v_i}+kw_1\p_{w_1}+kw_2\p_{w_2},\\
\bar\xi_j^1=-w_2\p_{v_j}+\sum_{i<j}v_{i-j+k}\p_{v_i}, \quad 1\leq j\leq k-1,\\
\bar\chi=-\bar\xi^2_1=-\sum_{i=1}^{k-2}(i+1)u_i\p_{u_i}+\sum_{i=1}^{k-1}(k-i)v_i\p_{v_i}
+kw_1\p_{w_1},\\
\bar\xi_j^2=\sum_{i<k-j}(i+j)u_{i+j-1}\p_{u_i}-\sum_{i<k-j+1}(k-i-j+1)v_{i+j-1}\p_{v_i},\quad 1<j\leq k-1,\\
\bar\xi^3_j=-(1-\delta_{j,1})w_2\p_{u_{k-j}}+\sum_{i<k-j}v_{i+j}\p_{u_i}+\delta_{j,1}w_2\p_{w_1},\quad 1\leq j\leq k-1.
\end{gather*}
Also let
\[
\bar\sigma=(\bar\xi_e+\bar\chi)/k=\sum_{i=1}^{k-1}v_i\p_{v_i}+w_2\p_{w_2}.
\]

We now test the above vector fields for tangency to $A=V(m^k_k)$.
Vector fields in $\Der(-\log D)$ preserve $\calF_1$, so for each $\xi\in\Der(-\log D)$ there exist $c_j\in\OO_T$, unique modulo $\mm_{T}$, such that
\[
\xi(m^k_k)=\sum_jc_jm^k_j.
\]
We determine their value modulo $\mm_T$ with the help of distinguished monomials. 
Let $\iota$ be the sign of the order-reversing permutation of $1,\ldots, k-1$.
Then, by Lemma \ref{28}, for $1<j\leq k$, the monomial $w_2^{k-2}v_{k-j+1}$ appears in the polynomial expansion of $m^k_j$ with coefficient $(-1)^{j-1}\iota$ but does not appear in the polynomial expansion of $m^k_\ell$ for $\ell\neq j$.
Similarly, the monomial $w_2^{k-1}$ has coefficient $\iota$ in the polynomial expansion of $m^k_1$ but does not appear in that of $m^k_j$ for $j\ge2$.

Let $\lambda_1',\dots,\lambda_k'$ denote the columns of the matrix $\lambda$ of \eqref{41} with its last row deleted. 
For any $\delta\in\Theta_T$, we have
\beq\label{117}
\delta(m^k_k)=\sum_{r=1}^{k-1}\det(\lambda_1',\dots,\delta(\lambda_r'),\dots,\lambda_{k-1}').
\eeq
For $\delta=\xi^2_j$, the only distinguished monomial to appear in any of the summands in \eqref{117} is
$v_jw_2^{k-2}$, which appears in the summand for $r=1$, with coefficient $(k-j)(-1)^k\iota$. 
Thus
\beq\label{120}
\xi^2_j(m^k_k)=(-1)^j(k-j)m^k_{k-j+1}\mod\mm_T\calF_1,\quad\text{for }1\leq j\leq k-1.
\eeq
Similarly we find
\[
\xi^1_j(\lambda_r')=
\begin{cases}
\lambda'_{k-j+r} & \text{if }r\leq j,\\
0 & \text{otherwise}.
\end{cases}
\]
Using \eqref{117} with $\delta=\xi^1_j$, it follows that
\beq\label{103}
\xi^1_j(m^k_i)=(-1)^{k-j+1}m^k_{i+j-k}\mod\mm_T\calF_1,\quad\text{for }1\le i\le k,\ 1\leq j\leq k-1.
\eeq
Note that \eqref{120}, and \eqref{103} with $i=k$, imply

\begin{prp}\label{104}
The map 
\beq\label{101}
dm^k_k:\Der(-\log D)\to \calF_1
\eeq
sending $\xi\in\Der(-\log D)$ to $\xi(m^k_k)$ is an $\calO_T$-linear surjection.
\qed
\end{prp}

Combining \eqref{120} and \eqref{103} with $i=k$, we construct vector fields
\[
\eta_j=(j-1)\xi^1_j-\xi^2_{k-j+1},\quad 2\leq j\leq k-1
\]
with linear part 
\[
\bar\eta_j\equiv(1-j)w_2\p_{v_j}+\sum_{i<j}(2j-i-1)v_{k+i-j}\p_{v_i}\mod\ideal{\p_u,\p_{w_1}},
\]
which lie in $\Der(-\log(D+A))$ to first order, since $\eta_j(m^k_k)\in
\mm_T\calF_1$. 

Both $\bar\chi$ and $\bar\sigma$ are semi-simple, and so by consideration of the distinguished monomials, $\chi$ and $\sigma$ must therefore lie in $\Der(-\log A)$ to first order.
The vector fields $\xi^3_j$ lie in $\Der(-\log A)$ to first order, since it is clear
by consideration of the distinguished monomials that $\xi^3_j(m^k_k)\in\mm_T\calF_1$. 

Thus we have $2k-1$ vector fields $\eta_2,\dots,\eta_{k-1},\chi,\sigma,\xi^3_1,\dots,\xi^3_{k-1}$ in $\Der(-\log D)$ which are also in $\Der(-\log A)$ to first order. 
By Proposition~\ref{104}, we can modify these by the addition of suitable linear combinations, with coefficients in $\mm_T$, of the Houston--Littlestone generators
of $\Der(-\log D)$, so that they are indeed in $\Der(-\log A)$ and therefore in $\Der(-\log (D+A))$. 
The determinant of the Saito matrix of the modified vector fields $\tilde\eta_2,\dots,\tilde\eta_{k-1},\tilde\chi,\tilde\sigma,\tilde\xi^3_1,\dots,\tilde\xi^3_{k-1}$ must be a multiple $\alpha hm^k_k$ of the equation of $D+A$. 
We now show that $\alpha$ is a unit, from which it follows, by Saito's criterion, that $D+A$ is a free divisor.

The modification of the vector fields does not affect the lowest order terms in the determinant of their Saito matrix, and these are the same as the lowest order terms in the determinant of the Saito matrix of their linear parts. 
With the rows representing the coefficients of $\p_{u_1},\dots,\p_{u_{k-2}},\p_{w_1},\p_{v_1},\dots,\p_{v_{k-1}},\p_{w_2}$ in this order, this matrix is of the form
\[
\begin{pmatrix}*&\bar B_1\\\bar B_2&0\end{pmatrix},
\]
with 
\[
\bar B_1=
\begin{pmatrix}
v_2&v_3&v_4&\cdots&v_{k-1}&-w_2\\
v_3&v_4&&\rdots&-w_2&0\\
v_4&&&\rdots&\rdots&\vdots\\
\vdots&\rdots&\rdots&&&\\
v_{k-1}&-w_2&\rdots&&&\vdots\\
w_2&0&\cdots&&\cdots&0
\end{pmatrix}
\]
and
\[
\bar B_2=
\begin{pmatrix}
2v_{k-1}&4v_{k-2}&6v_{k-3}&\cdots&\cdots&(2k-4)v_2&(k-1)v_1&v_1\\
-w_2&3v_{k-1}&5v_{k-2}&\cdots&\cdots&(2k-5)v_3&(k-2)v_2&v_2\\
0&-2w_2&4v_{k-1}&\ddots&&\vdots&\vdots&\vdots\\
\vdots&0&-3w_2&\ddots&\ddots&\vdots&\vdots&\vdots\\
&\vdots&0&\ddots&\ddots&kv_{k-2}&3v_{k-3}&v_{k-3}\\
\vdots&\vdots&\vdots&\ddots&\ddots&(k-1)v_{k-1}&2v_{k-2}&v_{k-2}\\
0&0&0&\cdots&0&-(k-2)w_2&v_{k-1}&v_{k-1}\\
0&0&0&\cdots&\cdots&0&0&w_2
\end{pmatrix}
\]
In its determinant we find the monomial $w_2^{2k-2}v_1$ with coefficient $\pm (k-1)!$. 
By Corollary~\ref{113}, this monomial is present in the equation of $D+A$. 
This proves that $\alpha$ is a unit and completes the proof that $D+V(m^k_k)$ is a free divisor.

The following consequence of Proposition~\ref{104} is needed to prove that Theorem~\ref{59} holds for \emph{any} adjoint divisor $A$ of $D$, and not just for $A=V(m^k_k)$.

\begin{cor}\label{111}
The adjoint divisor $A$ is unique up to isomorphism preserving $D$.
\end{cor}

\begin{proof} 
Let $A_0:=V(m^k_k)$. 
By Lemma~\ref{107}, any adjoint divisor $A_1$ must have an equation of the form $m_1:=m^k_k+\sum_{i=1}^{k-1}c_im^k_i$. 
Consider the family of divisors
\beq\label{112}
A:=V\left(m\right)\subset T\times(\CC,1),\quad m:=m_k^k+s\sum_{j=1}^{k-1}c_jm^k_j,
\eeq
where $s$ is a coordinate on $(\CC,1)$.
We claim that there exists a vector field 
\[
\Xi\in\Der_{T\times(\CC,1)/(\CC,1)}(-\log(D\times(\CC,1)))
\]
such that
\beq\label{119}
\Xi(m)=\p_s(m).
\eeq
Then the vector field $\p_s-\Xi$ is tangent to $D\times(\CC,1)$, and its integral flow trivializes the family \eqref{112}.
Passing to representatives of germs and scaling the $c_j$, the latter holds true in an open neighborhood of any point of the interval $\{0\}\times [0,1]\subset\CC^{n+1}\times\CC$.
A finite number of such neighborhoods cover this interval and it follows that $A_0$ and $A_1$ are isomorphic by an isomorphism preserving $D$.

To construct $\Xi$ we first show, using Proposition~\ref{104}, that
\beq\label{102}
dm_1:\Der(-\log D)\to\calF_1
\eeq
is surjective.
By \eqref{101}, there are vector fields $\delta_1,\dots,\delta_k\in\Der(-\log D)$, homogeneous with respect to the grading determined by the vector field $\chi$, such that $\delta_j(m^k_k)=m^k_j$ for $j=1,\dots,k$.
Pick $\alpha^\ell_{i,j}\in\calO_T$ such that 
\[
\delta_i(m^k_j)=\sum_{\ell=1}^k \alpha^\ell_{i,j}m^k_\ell,\quad 1\leq i,j\leq k,
\]
and set $L_j:=(\alpha^\ell_{i,j})_{1\leq\ell,i\leq k}$. 
The constant parts $L_j(0)$ are uniquely defined and $L_k$ is the identity matrix by choice of the $\delta_j$. 
For $j<k$, with respect to the grading determined by $\chi$, we have 
\[
\deg(m^k_k)<\deg(m^k_{k-1})\leq\cdots\leq\deg(m^k_1)\footnote{All the inequalities here are in fact strict, but we want to use the argument again later in a context where we assume only what is written here (see Corollary~\ref{115}).}.
\]
This gives
\[
\deg(\delta_i(m^k_j))>\deg(\delta_i(m^k_k))=\deg(m^k_i)\geq\deg(m^k_{i+1})\geq\cdots\geq\deg(m^k_k),
\]
and hence
\[
\delta_i(m^k_j)\in\ideal{m^k_{1},\ldots, m^k_{i-1}}+\mm_T\calF_1.
\] 
The constant matrices $L_j(0)$, $j<k$, are therefore strictly upper triangular, and the constant part $L_k+\sum_{j=1}^{k-1}c_j(0)L_j(0)$ of the matrix of $dm_1$ is invertible. 
Thus,
\[
dm_1(\Der(-\log D))+\mm_T\calF_1=\calF_1,
\]
and \eqref{102} follows by Nakayama's Lemma. 
As $m\equiv m_1\mod\mm_{T\times(\CC,1)}$, \eqref{102} gives
\[
dm(\Der_{T\times(\CC,1)/(\CC,1)}(-\log(D\times(\CC,1)))+\mm_{T\times(\CC,1)}\calO_{T\times(\CC,1)}\calF_1\supset\calO_{T\times(\CC,1)}\calF_1
\]
and Nakayama's Lemma yields
\beq\label{123}
dm(\Der_{T\times(\CC,1)/(\CC,1)}(-\log(D\times(\CC,1)))\supset\calO_{T\times(\CC,1)}\calF_1.
\eeq
Then any preimage $\Xi$ of $\p_s(m)\in\calF_1$ under $dm$ solves \eqref{119}.
\end{proof}

The proof of Theorem \ref{59} is now complete.
\end{proof}

\section{Discriminants of hypersurface singularities}\label{discr}

Let  $f\colon X:=(\CC^n,0)\to(\CC,0)=:T$ be weighted homogeneous of degree $d$ (with respect to positive weights) and have an isolated critical point at $0$.
Let $\chi_0$ be an Euler vector field for $f$, with $\chi_0(f)=d\cdot f$.
Denote by $J_f:=\ideal{\p_{x_1}(f),\dots,\p_{x_n}(f)}$ the Jacobian ideal of $f$.
Pick a weighted homogeneous $g=g_1,\dots,g_\mu\in\OO_X$ with decreasing degrees $d_i:=\deg(g_i)$ inducing a $\CC$-basis of the Jacobian algebra 
\[
M_f:=\calO_X/J_f.
\]
We may take $g_\mu:=1$ and $g_1:=H_0$ to be the Hessian determinant $H_0$ of $f$, which generates the socle of $M_f$.
Then 
\beq\label{121}
F(x,u):=f(x)+g_1(x)u_1+\cdots+g_\mu(x)u_\mu
\eeq
defines an $\calR_e$-versal unfolding 
\[
F\times\pi_S\colon Y:=X\times S\to T\times S
\]
of $f$, with base space $S:=(\CC^\mu,0)$, where
\[
\pi=\pi_S\colon Y=X\times S\to S
\]
is the natural projection.
Setting $\deg(u_i)=w_i:=d-d_i$ makes $F$ weighted homogeneous of degree $\deg(F)=d=\deg(f)$.
We denote by $\chi$ the Euler vector field $\chi_0+\delta_1$ where $\delta_1=\sum_{i=1}^\mu w_iu_i\p_{u_i}$.

Let $\Sigma\subset Y$ be the relative critical locus of $F$, defined by the relative Jacobian ideal $J_F^\rel=\ideal{\p_{x_1}(F),\dots,\p_{x_\mu}(F)}$, and set $\Sigma^0=\Sigma\cap V(F)$. 
Then $\calO_\Sigma$ is a finite free $\calO_S$-module with basis $g$.
As $\Sigma$ is smooth and hence Gorenstein, $\Hom_{\calO_S}(\calO_\Sigma,\calO_S)\cong\calO_\Sigma$ as $\calO_\Sigma$-modules, and a basis element $\Phi$ defines a symmetric perfect pairing 
\[
\ideal{\cdot,\cdot}\colon\calO_\Sigma\otimes_{\calO_S}\calO_\Sigma\to
\calO_S,\quad\ideal{g,h}:=\Phi(gh),
\]
which we refer to as the Gorenstein pairing.
As in \S\ref{images} (see the proof of Lemma~\ref{28}), a generator $\Phi$ may be defined by projection to the socle of the special fiber: 
We let $\Phi(h)$ be the coefficient of the Hessian $g_1$ in the expression of $h$ in the basis $g$.
By
\[
\ideal{\cdot,\cdot}_0\colon M_f\otimes_\CC M_f\to\CC
\]
we denote the induced Gorenstein pairing on $\calO_{\Sigma}/\mm_S\calO_{\Sigma}=\calO_X/J_f=M_f$.

Let $\check g=\check g_1,\dots,\check g_\mu$ denote the dual basis of $g$ with  respect to the Gorenstein pairing, and denote by $\check d_i$ the degree of $\check g_i$.
We have $d_i+\check d_i=d_1$, so $\check d_i=d_{\mu+i-1}$ (recall that we have ordered the $g_i$ by descending degree).

The discriminant $D=\pi_S(\Sigma^0)\subset S$ was shown by Kyoji Saito (see \cite{Sai80a}) to be a free divisor. 
The following argument proves this, and shows also that it is possible to choose 
a basis for $\Der(-\log D)$ whose Saito matrix is symmetric. 

\begin{thm}\label{1}
There is a free resolution of $\calO_{\Sigma^0}$ as $\calO_S$-module 
\[
\xymat{
0\ar[r] & \calO_S^\mu\ar[r]^\Lambda& \Theta_S\ar[r]^-{dF} & \calO_{\Sigma^0}\ar[r] & 0
}
\]
in which $\Lambda$ is symmetric, and is the Saito matrix of a basis of $\Der(-\log D)$.
\end{thm}

\begin{proof}
As in \eqref{tricky}, there is a commutative diagram with exact rows
\[
\xymat@C=10pt{
0\ar[r]&\calO_\Sigma\ar[rr]^-{F} && \calO_\Sigma\ar[rr] && \calO_{\Sigma^0}\ar[r] & 0\\
0\ar[r]&\calO_S^\mu\ar[dr]^-\cong\ar[u]^{\pp_{\check g}}_-\cong\ar[rr]^\Lambda &&\Theta_S \ar[u]^{\pp_g}_-\cong\ar[urr]_-{dF}\\
&&\Der(-\log D)\ar@{^{(}->}[ur]&
}
\]
where $\Lambda=(\lambda^i_j)_{1\le i,j\le\mu}$ is the matrix of multiplication by $F$ with respect to bases $\check g$ in the source and $g$ in the target. 
As in \eqref{106}, symmetry of $\Lambda$ follows from self-adjointness of multiplication by $F$ with respect to the Gorenstein pairing. 
Because of the form of $F$,  the map $\pp_g:\Theta_S\to\calO_{\Sigma}$ sending $\eta=\sum_j\alpha_j\p_{u_j}$ to $\sum_j\alpha_jg_j$ coincides with evaluation of $dF$ on a(ny) lift $\tilde\eta\in\Theta_Y$ of $\eta$; different lifts of the same vector field differ by 
a sum $\sum_j\alpha_j\p_{x_j}\in\Theta_{Y/S}$, and the evaluation of $dF$ on such a sum vanishes on $\Sigma$. 
The kernel of the composite $\Theta_S\to\OO_{\Sigma_0}$ consists of vector fields on $S$ which lift to vector fields on $Y$ which are tangent to $V(F)$, since $dF(\tilde \eta)$ is divisible by $F$ if and only if $\tilde\eta\in\Der(-\log V(F))$. 
It is well known (see e.g.\ \cite[Lem.~6.14]{Loo84}) that the set of  vector fields on $S$ which lift to vector fields tangent to $V(F)$ is equal to $\Der(-\log D)$.
\end{proof}

Denote by $m^i_j$ the $(\mu-1)$-minor of $\Lambda$ obtained by deleting the $i$th row and the $j$th column.
Then Lemmas~\ref{109}, \ref{110} and \ref{107} of \S\ref{images} remain valid in this new context.
That is, 
\begin{equation}\label{16}
\calF_1:=\calF_1^{\calO_S}(\calO_{\Sigma^0})=\ideal{m^\mu_j\mid j=1,\dots,\mu}_{\calO_S},\quad\calF_1\OO_{\Sigma^0}=\ideal{\mmm}_{\OO_{\Sigma^0}},
\end{equation}
and the adjoints of $D$ are divisors of the form
\[
A=V\bigl(\mmm+\sum_{j=1}^{\mu-1}c_jm^\mu_j\bigr).
\]

Although it is not part of the main thrust of our paper, the following result seems to be new, and is easily proved.
It assumes that $D$ is the discriminant of an $\calR_e$-versal deformation, but
does not require any assumption of weighted homogeneity. 
We denote by $H$ the relative Hessian determinant of the deformation \eqref{121}. 

\begin{thm}\label{12}
Let $A$ be any adjoint divisor for $D$. 
Then
\[
\tilde A:=\pi^{-1}(A)\cap\Sigma^0
\]
is a free divisor in $\Sigma^0$ containing $V(H)\cap \Sigma^0$, with reduced defining equation $(\mmm\circ\pi)/H$.
\end{thm}
 
\begin{proof}
By Corollary~\ref{115} below, we may assume that $A=V(\mmm)$, and hence $\tilde A=V(\mmm\circ\pi)\cap\Sigma_0$.
First, it is necessary to show that $H^2$ divides $\mmm\circ\pi$ and that $\mmm/H$ is reduced. 
Since $\Sigma^0$ is smooth, it is enough to check this at generic points of $V(H)$. 
This reduces to checking that it holds at an $A_2$-point. 
The miniversal deformation of an $A_2$-singularity is given by $G(x,v_1,v_2)=x^3+v_1x+v_2$.
In this case, $\mmm$ is, up to multiplication by a unit, simply the coefficient of $\p_{v_1}$ in the Euler vector field, namely $v_1$, and $v_1=-3x^2$ on $\Sigma^0$.
The Hessian $H$ is equal to $6x$, so $H^2$ does divide $\mmm\circ\pi$, and moreover the quotient $(\mmm\circ\pi)/H$ is reduced. 

As $\Sigma^0$ is the normalization of $D$ (see \cite[Thm.~4.7]{Loo84}), vector fields tangent to $D$ lift to vector fields tangent to $\Sigma^0$ (see \cite{Sei66}). 
Let $\tilde\delta_1,\dots,\tilde\delta_\mu$ be the lifts to $\Sigma^0$ of the symmetric basis $\delta_1,\dots,\delta_\mu$ of $\Der(-\log D)$ constructed
in Theorem~\ref{1}. 
We will show that $\tilde\delta_1,\dots,\tilde\delta_{\mu-1}$ form a basis for $\Der(-\log\tilde A)$. 

To see this, pick coordinates for $\Sigma^0$, and denote by $(\tilde\delta_1,\ldots,\tilde\delta_\mu)$ the matrix whose $j$th column consists of the coefficients of the vector field $\tilde\delta_j$ with respect to these coordinates.
Similarly, denote by $(\delta_1,\ldots,\delta_\mu)$ the matrix whose $j$th column consists of the coefficients of $\delta_j$ with respect to the coordinates $u_1,\ldots,u_\mu$.
We abbreviate $\pi_{\Sigma^0}:=\pi|_{\Sigma^0}$ and  $\pi_\Sigma:=\pi|_\Sigma$.
There is a matrix equality 
\beq\label{14}
[T\pi_{\Sigma^0}]\cdot(\tilde\delta_1,\dots,\tilde\delta_\mu)
=(\delta_1,\dots,\delta_\mu)\circ\pi,
\eeq
where $[T\pi_{\Sigma^0}]$ is the Jacobian matrix of $\pi$ with respect to the chosen coordinates. 
Let $\pi^\mu$ be obtained from $\pi$ by omitting the $\mu$th component, and let $(\delta_1,\dots,\delta_\mu)^\mu_\mu$ denote the submatrix of $(\delta_1,\dots,\delta_\mu)$ obtained by omitting its $\mu$th row and column. 
Then \eqref{14} gives
\[
[T\pi^{\mu}_{\Sigma^0}]\cdot(\tilde\delta_1,\dots,\tilde\delta_{\mu-1})=(\delta_1,\dots,\delta_\mu)^\mu_\mu\circ\pi,
\]
so that 
\beq\label{15}
\det[T\pi^\mu_{\Sigma^0}]\det(\tilde\delta_1,\ldots,\tilde\delta_{\mu-1})=m^{\mu}_\mu\circ\pi.
\eeq

We will now compute $\det[T\pi^\mu_{\Sigma^0}]$ in terms of $F$.
Because $g_\mu=1$, $u_\mu$ does not appear in the equations of $\Sigma$, so $\p_{u_\mu}\in T_{(x,u)}\Sigma$ for all $(x,u)\in\Sigma$. 
It follows that at any point of $\Sigma^0$, $T_{(x,u)}\Sigma$ has a basis consisting of a basis of $T_{(x,u)}\Sigma^0$ followed by the vector $\p_{u_{\mu}}$. 
With respect to such a basis, the matrix of $[T\pi_{\Sigma}]$ takes the form  
\[
[T\pi_{\Sigma}]=
\begin{bmatrix}
T\pi^{\mu}_{\Sigma^0}&*\\
0&1
\end{bmatrix}
\]
from which it follows that 
\beq\label{108}
\det[T\pi^\mu_{\Sigma^0}]=\det[T\pi_{\Sigma}].
\eeq
In order to express the latter in terms of $F$, we compare  two representations of the zero-dimensional Gorenstein ring 
\beq\label{999}
\OO_\Sigma/\ideal{u_1,\ldots,u_\mu}=\OO_\Sigma/\pi^*\mm_S=\OO_Y/\ideal{\p_{x_1}(F),\ldots,\p_{x_n}(F),u_1,\ldots,u_\mu}
\eeq
as a quotient of a regular $\CC$-algebra.
In both cases, by \cite[(4.7) Bsp.]{SS75}, the socle is generated by the Jacobian determinant of the generators of the respective defining ideal.
The first representation then shows that the socle is generated by $\det[T\pi_\Sigma]$, the second one, that it is also generated by the (relative) Hessian $H=\det(\p^2 F/\p x^2)$ of $F$. 
Hence, up to multiplication by a unit, we obtain 
\beq\label{114}
\det[T\pi_{\Sigma}]=H 
\eeq
in $\OO_{\Sigma}/\ideal{u_1,\ldots,u_\mu}$.
It is easy to see that the non-immersive locus of $\pi_{\Sigma}$ is the vanishing set of $H$. This, together
with \eqref{114}, shows that $\det[T\pi_{\Sigma^0}]=H$.
Now combining \eqref{15}, \eqref{108} and \eqref{114}, 
\[
\det(\tilde\delta_1,\dots,\tilde\delta_{\mu-1})=(\mmm\circ\pi)/H,
\]
and so is a reduced defining equation for $\tilde A=V(\mmm\circ\pi)$.
Finally, each $\tilde\delta_j$ is tangent to $\tilde A=(\pi_{\Sigma^0})^{-1}(A)$ at its smooth points, since $\delta_j$ is necessarily tangent to the non-normal locus $A\cap D$ of $D$.
The theorem now follows by Saito's criterion.
\end{proof}

\begin{rmk}\label{95}
Computation with examples appears to show that closure
\[
C_v:=\overline{\tilde A\setminus V(H)}
\]
is also a free divisor. 
\end{rmk}

We now go on to show first that the divisor $D+V(\mmm)$ is free and then (see Corollary~\ref{115}) that all adjoints are isomorphic.
Just as in \S\ref{images}, our proof makes use of the representation of $\Der(-\log D)$ on $\calF_1$, and relies on the surjectivity of $d\mmm\colon\Der(-\log D)\to \calF_1$.

\begin{prp}\label{13}
Assume that $d-d_1+2d_i\neq 0\neq d-d_i$ for $i=1,\dots,\mu$.
Then 
\[
d\mmm(\Der(-\log D))=\calF_1.
\]
\end{prp}

Inclusion of the left hand side in the right is a consequence of the $\Der(-\log D)$-invariance of $\calF_1$.
To show equality, it is enough to show that it holds modulo $\mm_S\calF_1$.
This will cover most of the remainder of this section.

Denote by $\bar \Lambda=(\bar \lambda^i_j)_{1\le i,j\le\mu}$ the linear part of $\Lambda$, and 
let $\bar\delta_i=\sum_j\bar \lambda^i_ju_i\p_{u_i}$ be the linear part of $\delta_i$.

\begin{thm}\label{3}
The entries of $\Lambda$ are given by $\lambda^i_j=\sum_{k=1}^\mu\ideal{\check g_i\check g_j,g_k}w_ku_k$.
In particular, $\bar \lambda^i_j=\sum_{k=1}^\mu\ideal{\check g_i\check g_j,g_k}_0w_ku_k$.
\end{thm}

\begin{proof}
Since $\chi_0(F)\in J_F^\rel$, we have
\[
F=\chi(F)\equiv\delta_1(F)=\sum_k w_ku_kg_k\mod J_F^\rel,
\]
and hence
\[
\lambda^i_j
=\ideal{\check g_i,F\check g_j}
=\ideal{\check g_i\check g_j,F}
=\ideal{\check g_i\check g_j,\sum_kw_ku_kg_k}
=\sum_k\ideal{\check g_i\check g_j,g_k}w_ku_k.
\]
\end{proof}

We call a homogeneous basis $g$ of $M_f$ self-dual if 
\beq\label{50}
\check g_i=g_{\mu+1-i}.
\eeq

\begin{lem}
$M_f$ admits self-dual bases.
\end{lem}

\begin{proof}
Denote by $W_j\subset M_f$ the subspace of degree-$d_j$ elements.
The space $W_1$ is $1$-dimensional generated by the Hessian of $f$.
Therefore $W_j$ and $W_k$ are orthogonal unless $d_j+d_k=d_1$, in which case $\ideal{\cdot,\cdot}_0$ induces a non-degenerate pairing $W_j\otimes_\CC W_k\to\CC$.
If $j\ne k$, one can choose the basis of $W_j$ to be the reverse dual basis of a basis of $W_k$.
Otherwise, $W_j=W_k$ and (since quadratic forms are diagonalizable) there is a basis of $W_j$ for which the matrix of $\ideal{\cdot,\cdot}_0$ is diagonal. 
Self-duality on $W_j$ is then achieved by a coordinate change with matrix 
\[
\begin{pmatrix}
1 & 0 & \cdots & \cdots & 0 & 1 \\
0 & \ddots & & & \rdots & 0 \\
\vdots & & 1 & 1 & & \vdots \\
\vdots & & i & -i & & \vdots \\
0 & \rdots & & & \ddots & 0\\
i & 0 & \cdots & \cdots & 0 & -i
\end{pmatrix}
\quad\text{or}\quad
\begin{pmatrix}
1 & 0 & \cdots & & \cdots & 0 & 1 \\
0 & \ddots & & & & \rdots & 0 \\
\vdots & & 1 & & 1 & & \vdots \\
 & & & 1 & & & \\
\vdots & & i & & -i & & \vdots \\
0 & \rdots & & & & \ddots & 0\\
i & 0 & \cdots & & \cdots & 0 &-i 
\end{pmatrix}
\] 
where $i=\sqrt{-1}$, for $\dim_\CC(W_j)$ even or odd, respectively.
A self-dual basis of $M_f$ is then obtained by joining the bases of the $W_j$ constructed above.
\end{proof}

Let $\bar m^i_j$ be the $(\mu-1)$-jet of $m^i_j$, that is, the corresponding minor of $\bar\Lambda$.

\begin{lem}\label{5}
Suppose $g$ is an $\calO_S$-basis for $\calO_\Sigma$ whose restriction to $M_f$ is self-dual, and that $d\ne d_i\ne0$. 
Then the following equalities hold true:
\begin{enumerate}[(a)]
\item\label{5a} $\mm_S\calF_1=\calF_1\cap\mm_S^\mu\calO_S$ and $\calF_1$ is minimally generated by $m^\mu_1,\dots,\mmm$.
In particular, $\bar m^\mu_i\equiv m^\mu_i\mod\mm_S\calF_1$ for $i=1,\dots,\mu$.
\item\label{5b} $\bar\delta_i(\bar m^\mu_\mu)=\pm(d-d_1+2d_i)\bar m^\mu_{\mu+1-i}$ for $i=2,\dots,\mu$.
\item\label{5c} $\delta_1(\mmm)\equiv\mmm\mod\CC^*$.
\end{enumerate}
\end{lem}

\begin{proof}
As $w_i=d-d_i\ne0$ by hypothesis, we may introduce new variables $v_i=w_iu_i$ for $i=1,\dots,\mu$. 
Under the self-duality hypothesis, Theorem~\ref{3} implies that the matrix $\bar \Lambda$ has the form
\beq\label{49}
\bar\Lambda=
\begin{pmatrix}
v_1 & v_2 & \cdots & \cdots & v_{\mu-1} & v_\mu \\
v_2 & \star & \cdots & \star & v_\mu & 0 \\
\vdots & \vdots & \rdots & \rdots & \rdots & \vdots \\
\vdots & \star & \rdots & \rdots & & \vdots \\
v_{\mu-1} & v_\mu & \rdots & & & \vdots \\
v_\mu & 0 & \cdots & \cdots & \cdots & 0
\end{pmatrix}
\eeq
where $\star$ entries do not involve the variable $v_\mu$. 
For the first row and column, this is clear. 
For the remaining entries, we note that, by Theorem~\ref{3}, $v_\mu$ appears in $\bar \lambda^i_j$ if and only if $0\ne\ideal{\check g_i\check g_j,g_\mu}_0=\ideal{\check g_i,\check g_j}_0$.
By the self-duality assumption \eqref{50}, this is equivalent to $i+j=\mu+1$, in which case $\ideal{\check g_i\check g_j,g_\mu}=1$ and $\bar \lambda_{\mu-i+1}^i=v_\mu$.

As in the proof of Theorem~\ref{59}, it is convenient to use $\iota$ to denote the sign of the order-reversing permutation of $1,\ldots,\mu-1$.
From \eqref{49} it follows that $\bar m^\mu_{\mu+1-i}$ involves a distinguished monomial $v_iv_\mu^{\mu-2}$, with coefficient $(-1)^{\mu-i-1}\iota$ for $i=1,\dots,\mu-1$ and $\iota$ for $i=\mu$; this monomial does not appear in any other minor $\bar m^\mu_{\mu+1-j}$ for $i\ne j$.
This implies \eqref{5a}.

In order to prove \eqref{5b}, assume, for simplicity of notation, that $\Lambda$ and $\delta$ are linear, and fix $i\in\{2,\dots,\mu-1\}$; the case $i=\mu$ is similar.
We know that $\delta_i(\mmm)$ is a linear combination of $m^\mu_1,\dots,\mmm$. 
We will show that
\beq\label{51}
\delta_i(\mmm)=(-1)^{i-1}(w_1-2w_{\mu-i+1})m^\mu_{\mu-i+1}
\eeq
by computing that the coefficient $c_{i,j}$ of the distinguished monomial $v_jv_\mu^{\mu-2}$ in $\delta_i(\mmm)$ satisfies
\beq\label{52}
c_{i,j}=(-1)^{\mu-2}\iota(w_1-2w_{\mu-i+1})\delta_{i,j}.
\eeq
The self-duality assumption \eqref{50} implies that $d_1-d_\ell=\check d_\ell=d_{\mu-\ell+1}$.
Using $w_\ell=d-d_\ell$, this gives
\[
w_1-2w_{\mu-i+1}=d-d_1+2d_{\mu-i+1}-2d=d-d_1+2d_1-2d_i-2d=-(d-d_1+2d_i).
\]
So \eqref{5b} will follow from \eqref{51}.

By linearity of $\delta_i$, the only monomials in the expansion of $\mmm$ that could conceivably contribute to a non-zero $c_{i,j}$ are of the following three forms:
\beq\label{53}
v_\mu^{\mu-1},\quad 
v_jv_\mu^{\mu-2},\quad 
v_jv_kv_\mu^{\mu-3}.
\eeq
The first monomial does not figure in the expansion of $\mmm$. 
Monomials of the other two types do appear.
The second type of monomial in \eqref{53} must satisfy $j=1$ and arises as the product
\beq
\label{54}
(-1)^{\mu-2}\iota v_1v_\mu^{\mu-2}=(-1)^{\mu-2}\iota \lambda^1_1\lambda^2_{\mu-1}\lambda^3_{\mu-2}\cdots \lambda^{\mu-1}_2.
\eeq
Monomials of the third type in \eqref{53} must satisfy $k=\mu-j+1$. 
Each such monomial arises in the expansion of $\mmm$ in two ways, which coincide when $j=\mu-j+1$:
\begin{align}
(-1)^{\mu-3}\iota v_jv_kv_\mu^{\mu-3}&=(-1)^{\mu-3}\iota\lambda^1_j\lambda^2_{\mu-1}\cdots \lambda^{\mu-j}_{j+1}\lambda^{\mu-j+1}_1\lambda^{\mu-j+2}_{j-1}\cdots \lambda^{\mu-1}_2\label{55}\\
(-1)^{\mu-3}\iota v_jv_kv_\mu^{\mu-3}&=(-1)^{\mu-3}\iota\lambda^1_{\mu-j+1}\lambda^2_{\mu-1}\cdots\lambda^{j-1}_{\mu-j+2}\lambda^j_1\lambda^{j+1}_{\mu-j}\cdots \lambda^{\mu-1}_2\label{56}.
\end{align}

In terms of the coordinates $v_1,\dots,v_\mu$, $\delta_i$ contains monomials
\begin{align}
w_iu_i\p_{u_1}&=w_1v_i\p_{v_1}\label{57},\\
w_\mu u_\mu\p_{u_{\mu-i+1}}&=w_{\mu-i+1}v_\mu\p_{v_{\mu-i+1}}\label{58}.
\end{align}
Now \eqref{57} applied to \eqref{54} contributes $w_1(-1)^{\mu-2}\iota$ to $c_{i,i}$, \eqref{58} applied to one or two copies of \eqref{55} for $i=j$ contributes $2(-1)^{\mu-3}\iota w_{\mu-i+1}$ to $c_{i,i}$ in both cases.
There are no contributions to the coefficient of any other distinguished monomial. 

We have proved \eqref{51}, from which \eqref{5b} follows; \eqref{5c} is clear, since $\delta_1$ is the Euler vector field.
\end{proof}  

By Nakayama's lemma, Proposition~\ref{13} is now an immediate consequence of \eqref{16} and Lemma~\ref{5}.

The next result, closely analogous to Corollary~\ref{111}, follows from Proposition~\ref{13} by the same argument by which Corollary~\ref{111} is deduced from Proposition~\ref{104}.  

\begin{cor}\label{115}
Assume the hypothesis of Proposition~\ref{13}. 
Then any two adjoint divisors of $D$ are isomorphic by an isomorphism preserving $D$.\qed
\end{cor}

\begin{prp}\label{106}
Let $D=V(h)$ and $A=V(m)$ be divisors in $S$, and suppose that $D$ is a free divisor. 
Let $\calF$ be the $\OO_S$-ideal $dm(\Der(-\log D))$, and suppose that $m\in\calF$. 
Then the following two statements are equivalent:
\begin{enumerate}
\item $\depth_{\OO_S}\calF=\mu-1$. 
\item $D+A$ is a free divisor.
\end{enumerate}
\end{prp}

\begin{proof}\pushQED{\qed}
Apply the depth lemma (see \cite[Prop.~1.2.9]{BH93}) to the two short exact sequences:
\begin{gather*}
\xymat{0\ar[r]&\ideal{m}\ar[r]&\calF\ar[r]&\calF/\ideal{m}\ar[r]&0}\\
\xymat{0\ar[r]&\Der(-\log (D+A))\ar[r]&\Der(-\log D)\ar[r]^-{dm}&\calF/\ideal{m}\ar[r]&0}\qedhere
\end{gather*}
\end{proof}

\begin{proof}[Proof of Theorem~\ref{60}]
By Corollary~\ref{115}, we may assume that $A=V(\mmm)$.
Recall that $d\mmm(\Der(-\log D))=\calF_1$ by Proposition~\ref{13} and hence $\depth_{\OO_S}\calF_1=\mu-1$ by the Hilbert--Burch Theorem (see \cite[Thm.~1.4.17]{BH93}).
So Proposition~\ref{106} with $m=\mmm$ and $\calF=\calF_1$ yields the claim.
\end{proof}

\begin{rmk}
We remark that the very simple deduction of Theorem \ref{60} from Propositions \ref{13} and \ref{106}
does not have a straightforward analogue by which Theorem \ref{59} can be deduced from Propositions \ref{104} and \ref{106}. 
Firstly, the image $D$ of a stable map-germ $f:(\CC^n,0)\to(\CC^{n+1},0)$ is not free: $\Der(-\log D)$ has depth $n$ and not $n+1$. 
Secondly, there can be no way of adapting the argument to deal with this difference without some other input, since when $f$ is the stable germ of corank $2$ of Example \ref{17b}, the corresponding map $\Der(-\log D)\to\calF_1$ is surjective, and $\calF_1$ has depth $n$, but even so $\Der(-\log (D+A))$ is not free.
\end{rmk}

We conclude this section with a description of the relation between the adjoint divisor of $D$ and the bifurcation set of the deformation.
For $u\in S$, we set $X_u:=\pi_S^{-1}(u)$ and define $f_u\colon X_u\to T$ by $f_u(x):=F(x,u)$.
We consider $S':=(\CC^{\mu-1},0)$ with coordinates $u'=u_1,\ldots,u_{\mu-1}$, and we denote by
\begin{equation}\label{90}
\rho\colon S\to S'\quad u\mapsto u',
\end{equation}
the natural projection forgetting the last coordinate.
Recall that the \emph{bifurcation set} is the set $B\subset S'$ of parameter values $u'$ such that $f_{u'}:=f_{(u',0)}$ has fewer than $\mu$ distinct critical values. 
The coefficient $u_\mu$ of $g_\mu=1$ is set to $0$ since it has no bearing on the number of critical values. 
The bifurcation set consists of two parts: the level bifurcation set $B_v$ consisting of parameter values $u'$ for which $f_{u'}$ has distinct critical points with the
same critical value, and the local bifurcation set $B_\ell$ where $f_{u'}$ has a degenerate critical point.
H.~Terao, in \cite{Ter83}, and J.W.~Bruce in \cite{Bru85} proved that $B$ is a free divisor and gave algorithms for constructing a basis for $\Der(-\log B)$.
The free divisor $B$ is of course singular in codimension $1$. 
The topological double points (points at which $B$ is reducible) are of four generic types:
\begin{asparaitem} 
\item Type 1: $f_{u'}$ has two distinct degenerate critical points, $x_1$ and $x_2$.
\item Type 2: $f_{u'}$ has two distinct pairs of critical points, $x_1$, $x_2$ and $x_3,x_4$, such that $f_{u'}(x_1)=f_{u'}(x_2)$ and $f_{u'}(x_3)=f_{u'}(x_4)$.
\item Type 3: $f_{u'}$ has a pair of critical points $x_1$ and $x_2$ with the same critical value, and also a degenerate critical point $x_3$.
\item Type 4: $f_{u'}$ has three critical points $x_1,x_2$ and $x_3$ with the same critical value.
\end{asparaitem}
In the neighborhood of a double point of type 1, 2 or 3, $B$ is a normal crossing of two smooth sheets. 
In the neighborhood of a double point of type 4, $B$ is isomorphic to a product $B_0\times(\CC^{\mu-2},0)$ where $B_0=V(uv(u-v))\subset(\CC^2,0)$. 

\begin{prp} 
For any adjoint divisor $A$ for $D$, \eqref{90} induces a surjection
\begin{equation}\label{91}
\rho\colon D\cap A\onto B.
\end{equation}
\end{prp}

\begin{proof}
A point $u\in S$ lies in $D\cap A$ if the sum of the lengths of the Jacobian algebras of $f_u$ at  points $x\in f_u^{-1}(0)$ is greater than $1$. 
The sum may be greater than $1$ because for some $x$ the dimension of the Jacobian algebra is greater than 1 -- in which case $f_u$ has a degenerate critical point at $x$ -- or because $f_u$ has two or more critical points with critical value $0$. 
In either case, it is clear that $\rho(u)\in B$. 
If $u'\in B$, then $f_{u'}$ has either a degenerate critical point or a repeated critical value (or both). 
In both, cases let $v$ be the corresponding critical value.
Then $(u',-v)\in D\cap A$, proving surjectivity.
\end{proof}

\begin{rmk}
The projection \eqref{91} is a partial normalization, in the sense that generically, topological double points of $u'\in B$ of types 1, 2 and 3 are separated.
Indeed, in each such case $f_{u'}$ has two critical points with different critical values, and hence with different preimages under $\rho$.
However, a general point $u'$ of type 4 has only one preimage, $(u',-f_{(u',0)}(x_i))$, in $D\cap A$. 
Generically, at such a point $D$ is a normal crossing of three smooth divisors, and $D\cap A$ is the union of their pairwise intersections. 
Thus $B_0$ is improved to a curve isomorphic to the union of three coordinate axes in $3$-space. 
\end{rmk}

Finally, our free divisors $D+A$ and $\tilde A$ of Theorems~\ref{60} and \ref{12}, and the conjecturally free divisor $C_v$ of Remark~\ref{95}, fit into the following commutative diagram, in which $A$ is any adjoint divisor for $D$, and the simple and double underlinings indicate conjecturally free and free divisors.
\[
\xymat{
V(H)\cap\Sigma^0\ar@/^1pc/@{^(->}[rr]\ar@{>>}[dd]_-{\rho\circ\pi|}&\underline{C_v}\ar@{^(->}[r]\ar@{>>}[dd]_{\rho\circ\pi|}&\underline{\underline{\tilde A}}\ar[d]^-\pi\ar@{^(->}[r]&\Sigma^0\ar[d]^-\pi\\
&&D\cap A\ar@{>>}[d]_-\rho\ar@{^(->}[r]&\underline{\underline{D}}\ar@{^(->}[r]&\underline{\underline{D+A}}\\
B_\ell\ar@/_1pc/@{^(->}[rr]&B_v\ar@{^(->}[r]&\underline{\underline{B}}
}
\]

\section{Pull-back of free divisors}\label{lfd}

In this section, we describe a procedure for constructing new free divisors from old by a pull-back construction. 
It is motivated by Example \ref{17}.(\ref{17d}). 

\begin{thm}\label{19}
Suppose that $D=\bigcup_{i=1}^kD_i\subset(\CC^n,0)=:X$ is a germ of a free divisor.
Let $f\colon X\to Y:=(\CC^k,0)$ be the map whose $i$th component $f_i\in\OO_X$, for $i=1,\dots,k$, is a reduced equation for $D_i$. 
Suppose that, for $j=1,\dots,k$, there exist vector fields $\eps_j\in\Theta_X$ such that
\beq\label{127}
df_i(\eps_j)=\delta_{i,j}\cdot f_i.
\eeq
Let $N:=V(y_1\cdots y_k)\subset Y$ be the normal crossing divisor, so that $D=f^{-1}(N)$. 
Let $E\subset Y$ be a divisor such that $N+E$ is free.
Then provided that no component of $f^{-1}(E)$ lies in $D$, $f^{-1}(N+E)=D+f^{-1}(E)$ is a free divisor.
\end{thm}

\begin{proof}
The vector fields $\eps_1,\dots,\eps_k$ can be incorporated into a basis $\eps_1,\dots,\eps_n$ for $\Der(-\log D)$ such that \eqref{127} holds for $j=1,\dots,n$ and hence
\beq\label{18}
tf(\eps_j)=\sum_{i=1}^kdf_i(\eps_j)\p_{y_i}=
\begin{cases}
\omega f(y_j\p_{y_j}), & \text{ if }j\le k,\\
0, & \text{ otherwise}.
\end{cases}
\eeq

Any Saito matrix of $N+E$ can be written in the form $S_{N+E}=S_N\cdot A$, where $S_N=\diag(y_1,\dots,y_k)$ is the standard Saito matrix of $N$ and $A=(a_{i,j})\in\calO_Y^{k\times k}$. 
Then, by Saito's criterion, $h:=\det A$ and $g:=y_1\cdots y_nh$ are reduced equations for $E$ and $N+E$ respectively.
For $j=1,\dots k$, consider the vector fields $v_j:=\sum_{i=1}^ka_{i,j}y_i\p_{y_i}\in\Der(-\log(N+E))$ whose coefficients are the columns of $S_{N+E}$ and let \[
\tilde{v}_j:=\sum_{i=1}^k(a_{i,j}\circ f)\eps_i\in\Theta_X.
\] 
By \eqref{18}, we have $tf(\tilde{v}_j)=\omega f(v_j)$; by construction of the $v_j$, it follows that 
\[
d(g\circ f)(\tilde v_j)=(dg(v_j))\circ f\in (g\circ f)\calO_X,
\]
so that $\tilde v_j\in\Der(-\log f^{-1}(N+E))$, and moreover $\eps_j\in\Der(-\log f^{-1}(N+E))$ for $j>k$.

Let $S_D$ be the Saito matrix of $D$, whose columns are the coefficients of the vector fields $\eps_1,\dots, \eps_n$. 
The matrix of coefficients of the vector fields $\tilde v_1,\dots,\tilde v_k,\eps_{k+1},\dots,\eps_n$ is equal to
\[
S_D
\cdot
\begin{pmatrix}A\circ f&0\\
0&I_{n-k}
\end{pmatrix}
\]
and thus its determinant $\det(S_D)\cdot(h\circ f)$ defines $D+f^{-1}(E)=f^{-1}(N+E)$. 
By Saito's criterion, this shows that the latter is a free divisor, provided 
\ben
\y
$h\circ f$ is reduced and 
\y
$h\circ f$ has no irreducible factor in common with $f_1\cd f_k$. 
\een
Now $h\circ f$ is reduced where $f$ is a submersion. 
So provided no component of
$f^{-1}(E)$ is contained in the critical set $\Sigma_f$ of $f$, $h\circ f$ is reduced.
In fact $\Sigma_f=D_{\text{\scriptsize Sing}}$.
To see this, consider the ``logarithmic Jacobian matrix'' $(\eps_i(f_j))_{1\leq i\leq n, 1\leq j\leq k}$ of $f$.
The determinant of the first $k$ columns is equal to $f_1\cdots f_k$; thus $f_1\cdots f_k$ is in the Jacobian ideal of $f$, so $\Sigma_f\subset D$. Thus $D_{\text{\scriptsize Sing}}=D\cap\Sigma_f=\Sigma_f$. This shows that condition (2) above implies  condition (1).
\end{proof}

\begin{exa}
Let $D:=V(x_1\cdots x_n)\subset(\CC^n,0)$ be the normal crossing divisor, define $f\colon(\CC^n,0)\to(\CC^2,0)$ by $f(x_1,\ldots,x_n)=(x_1\cdots x_k, x_{k+1}\cdots x_n)$, and let $g\in\OO_{\CC^2,0}$ be any germ not divisible by either of the coordinates. 
Take $N:=V(y_1y_2)$ and $E:=V(g(y_1,y_2))$. 
By Theorem~\ref{19} it follows that $V(x_1\cd x_ng(x_1\cd x_k,x_{k+1}\cd x_n))$ is a free divisor. 
The condition that no component of $f^{-1}(E)$ should lie in $D$ is guaranteed by the requirement that neither $y_1$ nor $y_2$ should divide $g(y_1,y_2)$. 
\end{exa}

The existence hypothesis \eqref{127} in the theorem is not fulfilled for every reducible free divisor. 
In the graded case, the vector fields $\eps_j$ must have degree zero.  
If $D$  is the discriminant of a versal deformation of a weighted homogeneous isolated hypersurface singularity meeting the hypotheses of Theorem \ref{60}, and $A$ is an adjoint divisor, then $D+A$ is free but the only vector fields of weight zero in $\Der(-\log D)$ are multiples of the Euler field. 
Thus, hypothesis \eqref{127} cannot hold. 

However there is an interesting class, namely linear free divisors, for which this requirement always holds. 
We recall from \cite{GMNS09} that a free divisor $D$ in the $n$-dimensional vector space $V$ is said to be {linear} if $\Der(-\log D)$ has a basis consisting of vector fields of weight $0$. 
The linear span $\Der(-\log D)_0$ of the basic fields is an $n$-dimensional Lie algebra. 
It is naturally identified with the Lie algebra of the algebraic subgroup $\iota\colon G_D\into\GL(V)$ consisting of the identity component of the set of automorphisms preserving $D$. 
It follows that $(V,G_D,\iota)$ is a prehomogeneous vector space (see \cite{SK77}) with discriminant $D$. 

Let $D\subset V$ be a linear free divisor and $D=\bigcup_{i=1}^kD_i$ a decomposition into irreducible components.
The corresponding defining equations $f_1,\dots,f_k$ are polynomial relative invariants of $(V,G_D,\iota)$ with associated characters $\chi_1,\dots,\chi_k$; that is, for $g\in G_D$ and $x\in V$, $f_j(gx)=\chi_j(g)f_j(x)$. These characters are multiplicatively independent, by \cite[\S 4, proof of Prop.~5]{SK77}.
Let $\fg_D$ denote the Lie algebra of $G_D$.  
By differentiating the character map $\chi=(\chi_1,\dots,\chi_k)\colon G\to(\CC^*)^k$ we obtain an epimorphism of Lie algebras $d\chi\colon\fg_D\to\CC^k$.
This yields a decomposition
\[
\fg_D=\ker d\chi\oplus\bigoplus_{i=1}^k\CC\eps_i,\quad d\chi_i(\eps_j)=
\delta_{i,j}.
\]
For $\delta\in\fg_D$, the equality $f_i(gx)=\chi_i(g)f_i(x)$ differentiates to $\delta(f_i)=d\chi_i(\delta)\cdot f_i$, which implies \eqref{127}.
We proved

\begin{prp}\label{20}
Any germ of a linear free divisor $D$ satisfies the existence hypothesis \eqref{127} of Theorem~\ref{19}.\qed
\end{prp}

\begin{exa}
Let $\sigma_i(y)$ be the $i$th symmetric function of $y=y_1,\dots,y_k$ and set $N:=V(\sigma_k(y))$ and $E:=V(\sigma_{k-1}(y))$.
As seen in Example~\ref{17}.\eqref{17d}, the divisor $N+E$ is free.  
So by Theorem~\ref{19} and Proposition~\ref{20}, for any germ $D$ of a linear free divisor with distinct irreducible components $D_i=V(f_i)$, the divisor germ $V(\sigma_k(f_1,\ldots, f_k)\sigma_{k-1}(f_1,\ldots,f_k))$ is also free. 
No component of $V(\sigma_{k-1}(f_1,\ld,f_k))$ can lie in $D$, since were this the case, some $f_i$ would divide $\sigma_{k-1}(f_1,\ldots,f_k)$ and therefore would divide $f_1\cdots\widehat{f_i}\cdots f_k$.  

If each of the $D_i$ is normal, then in fact $f^{-1}(E)$ is an adjoint divisor of the normalization of $D$.  
As the singular locus of any free divisor has pure codimension $1$, the singular locus of $D$ is equal to its non-normal locus.
The ring of functions on the normalization $\bar D=\coprod_{i=1}^kD_i$ has presentation matrix $\diag(f_1,\dots,f_k)$. 
Thus $V(\sigma_{k-1}\circ f)$ is an adjoint divisor of $D$. 
\end{exa}

There is another class of divisors that fits naturally into the setup of Theorem~\ref{19}, namely that of hyperplane arrangements.

\begin{prp}\label{21}
Given the hypothesis \eqref{127}, any germ $N+E$ of a free hyperplane arrangement automatically satisfies the hypothesis on $f^{-1}(E)$ in Theorem~\ref{19}.
\end{prp}

\begin{proof}
By assumption, $N=V(y_1,\dots,y_k)$ and $E=\bigcup_{i=k+1}^mH_i$ where $H_i=V(\ell_i)$ for some linear equation $\ell_i(y)=\sum_j\alpha_{i,j}y_j$ for $i=k+1,\dots,m$.
We need to show that no component of any of the $\ell_i\circ f$ is divisible by any $f_j$. 
Suppose to the contrary that $\ell_i\circ f=g\cdot f_t$.
For $s\ne t$, $\eps_s$ applied to this equation gives $\alpha_{i,s}f_s=\eps_s(g)\cdot f_t$.
Since $f_t$ cannot divide $f_s$ as $D$ is reduced, it follows that $\alpha_{i,s}=0$ for any $s\neq t$, and thus $\ell_i=\alpha_{i,t}y_t$. 
This is absurd since $A$ is supposed reduced. 
\end{proof}

Combining Propositions~\ref{20} and \ref{21} proves

\begin{cor}
Let $A=\bigcup_{i=1}^mH_i\subset(\CC^k,0)$ be the germ of a free hyperplane arrangement containing the normal crossing divisor $\{y_1\cd y_k=0\}$,
and let $D\subset(\CC^n,0)$ be the germ of a linear free divisor whose irreducible components have equations $f_1,\ldots,f_k$. 
If $f\colon(\CC^n,0)\to(\CC^k,0)$ is defined by $f(x)=(f_1(x),\ldots,f_k(x))$, then
$f^{-1}(A)$ is a free divisor.
\end{cor} 

Note that the corollary applies to any essential free arrangement, once suitable coordinates are chosen.

\bibliographystyle{amsalpha}
\bibliography{fad}
\end{document}